\documentclass[11pt]{article}
\usepackage{amsmath}
\usepackage{amssymb}
\usepackage{amsfonts}
\usepackage{eucal}
\usepackage{graphicx,url}
\usepackage{amssymb,amsmath, amsthm}
\usepackage{epstopdf}
\newtheorem{Theorem}{Theorem}[section]
\newtheorem{Proposition}[Theorem]{Proposition}

\begin{document}
\title{Large deviations for fractional Poisson processes\thanks{The
financial support of the Research Grant PRIN 2008
\emph{Probability and Finance} is gratefully acknowledged.}}
\author{Luisa Beghin\thanks{Address: Dipartimento di
Scienze Statistiche, Sapienza Universit\`a di Roma, Piazzale Aldo
Moro 5, I-00185 Roma, Italy. e-mail:
\texttt{luisa.beghin@uniroma1.it}}\and Claudio
Macci\thanks{Dipartimento di Matematica, Universit\`a di Roma Tor
Vergata, Via della Ricerca Scientifica, I-00133 Roma, Italia.
E-mail: \texttt{macci@mat.uniroma2.it}}}
\date{}
\maketitle
\begin{abstract}
\noindent We prove large deviation principles for two versions of
fractional Poisson processes. Firstly we consider the \emph{main
version} which is a renewal process; we also present large
deviation estimates for the ruin probabilities of an insurance
model with constant premium rate, i.i.d. light tail claim sizes,
and a fractional Poisson claim number process. We conclude with
the \emph{alternative version} where all the random variables are
weighted Poisson distributed.\\


\end{abstract}

\section{Introduction}
The theory of large deviations gives an asymptotic computation of
small probabilities on exponential scale (we refer to
\cite{DemboZeitouni} for this topic). The aim of this paper is to
prove some large deviation results for fractional Poisson
processes.  To the best of our knowledge, these techniques have
not been applied so far to the fractional Poisson process.

The study of fractional versions of the usual renewal processes
has recently received an increasing interest, starting from the
paper by \cite{RepinSaichev}. In \cite{Jumarie} the so-called
fractal Poisson process is introduced (by means of non-standard
analysis) and proved to have independent increments. In analogy
with the fractional Brownian motion, in
\cite{WangWen}-\cite{WangWenZhang} is proposed a process
constructed as a stochastic integral with respect to the Poisson
measure. A different approach has been followed by other authors,
in the mainstream of the fractional diffusions, in the sense of
extending well-known differential equations by introducing
fractional-order derivatives with respect to time: the relaxation
equation (see e.g. \cite{Nonnenmacher}), the heat and wave
equations (see e.g.
\cite{Fujita}-\cite{Mainardi-AML-1996}-\cite{Mainardi-CSF-1996})
as well as the telegraph equation (see e.g.
\cite{OrsingherBeghin2004}) and the higher-order heat-type
equations (see e.g. \cite{Beghin}). In this context the solution
of the fractional generalization of the Kolmogorov-Feller
equation, together with the distribution of the waiting time for
the corresponding process, is derived in \cite{Laskin2003}. Many
other aspects of this type of fractional Poisson processes have
been analyzed: a probabilistic representation of the fractional
Poisson process of order $\nu$ as a composition of a standard
Poisson process with a random time given by a fractional diffusion
is given in \cite{BeghinOrsingher2009} (note that this has some
analogy with other results holding for compositions of different
processes; see e.g.
\cite{BeghinOrsingherSakhno}-\cite{OrsingherBeghin2009}) and, for
$\nu=1/2$, the time argument reduces to the absolute value of a
Brownian motion; in \cite{MeerschaertNaneVellaisamy} it is proved
that we have the same one-dimensional distributions of a standard
Poisson process time-changed via an inverse stable subordinator;
in \cite{PolitiKaizojiScalas} it is given a full characterization
of the fractional Poisson process in terms of its multidimensional
distributions. Other aspects of the fractional Poisson process are
analyzed in
\cite{MainardiGorenfloScalas}-\cite{MainardiGorenfloVivoli}-\cite{UchaikinCahoySibatov}.

We also recall other references with different approaches.
Applications based on fractional Poisson processes can be found in
\cite{UchaikinSibatov} (in the field of the transport of charge
carriers in semiconductors) and in \cite{Laskin2009} (in the field
the fractional quantum mechanics); an inference problem is studied
in \cite{CahoyUchaikinWoyczynski}; a version of space-fractional
Poisson process where the state probabilities are governed by
equations with a fractional difference operator found in time
series analysis is presented in \cite{OrsingherPolito}.

The outline of the paper is the following. We start with some
preliminaries in section \ref{sec:preliminaries}. In section
\ref{sec:main-version} we consider the \emph{main version} which
is a slight generalization of the renewal process in
\cite{BeghinOrsingher2009}-\cite{BeghinOrsingher2010}. We give
results for the empirical means of the i.i.d. holding times and
for the normalized counting processes; furthermore we study an
insurance model with fractional Poisson claim number process. In
section \ref{sec:alternative-version} we present large deviation
results for an \emph{alternative version}, which is the first
version presented in section 4 in \cite{BeghinOrsingher2009} with
a suitable deterministic time-change. In such a case we have a
weighted Poisson process, i.e. all the random variables are
weighted Poisson distributed with the same weights (namely the
weights do not depend on $t$). In the literature a weighted
Poisson process is defined in \cite{BalakrishnanKozubowski} and
examples of weighted Poisson distributions can be found in
\cite{CastilloPerez1998} and \cite{CastilloPerez2005}.

\section{Preliminaries}\label{sec:preliminaries}

\paragraph{Preliminaries on large deviations.} We start by
recalling some basic definitions (see \cite{DemboZeitouni}, pages
4-5). Given a topological space $\mathcal{Z}$, we say that a
family of $\mathcal{Z}$-valued random variables $\{Z_t:t>0\}$
satisfies the large deviation principle (LDP from now on) with
rate function $I$ if: the function $I:\mathcal{Z}\to [0,\infty]$
is lower semi-continuous; the upper bound
$$\limsup_{t\to\infty}\frac{1}{t}\log P(Z_t\in C)\leq-\inf_{x\in C}I(x)$$
holds for all closed sets $C$; the lower bound
$$\liminf_{t\to\infty}\frac{1}{t}\log P(Z_t\in G)\geq-\inf_{x\in G}I(x)$$
holds for all open sets $G$. The above definition can be given
also for a sequence of $\mathcal{Z}$-valued random variables
$\{Z_n:n\geq 1\}$ (we mean the discrete parameter denoted by $n$
in place of the continuous parameter $t$). Moreover a rate
function is said to be good if all its level sets
$\{\{x\in\mathcal{Z}:I(x)\leq\eta\}:\eta\geq 0\}$ are compact.
Throughout this paper we often refer to the well-known G\"{a}rtner
Ellis Theorem (see e.g. Theorem 2.3.6 in \cite{DemboZeitouni}).
Furthermore we always set $0\log 0=0$ and $0\log\frac{0}{0}=0$.

\paragraph{Preliminaries on (generalized) Mittag Leffler function.}
The Mittag Leffler function is defined by
$$E_{\alpha,\beta}(x):=\sum_{r\geq 0}\frac{x^r}{\Gamma(\alpha r+\beta)}$$
(see e.g. \cite{Podlubny}, page 17). We recall that, if we write
$a_t\sim b_t$ to mean that $\frac{a_t}{b_t}\to 1$ as $t\to\infty$,
we have
\begin{equation}\label{eq:asymptotic}
E_{\nu,\beta}(z)\sim\frac{1}{\nu}z^{(1-\beta)/\nu}e^{z^{1/\nu}}\
\mathrm{as}\ z\to\infty
\end{equation}
(see e.g. eq. (1.8.27) in \cite{KilbasSrivastavaTrujillo}).
Finally we recall that the generalized Mittag Leffler function is
defined by
$$E_{\alpha,\beta}^\gamma(x):=\sum_{r\geq 0}\frac{(\gamma)_rx^r}{r!\Gamma(\alpha r+\beta)},$$
where $(\gamma)_r=1$ is the Pochammer symbol defined by
$$(\gamma)_r:=\left\{\begin{array}{ll}
1&\ \mathrm{if}\ r=0\\
\gamma(\gamma+1)\cdots(\gamma+r-1)&\ \mathrm{if}\
r\in\{1,2,3,\ldots\}
\end{array}\right.$$
(see e.g. eq. (1.9.1) in \cite{KilbasSrivastavaTrujillo}); note
that $E_{\alpha,\beta}^1$ coincides with $E_{\alpha,\beta}$.

\section{Results for the main version (renewal process)}\label{sec:main-version}
Throughout this section we consider a class of fractional Poisson
processes defined as renewal processes. More precisely, for
$\nu\in(0,1]$ and $h,\lambda>0$, we consider
$\{M_{\nu,h,\lambda}(t):t\geq 0\}$ defined by
\begin{equation}\label{eq:def-main-version}
M_{\nu,h,\lambda}(t):=\sum_{n\geq 1}1_{T_1+\cdots+T_n\leq t},
\end{equation}
where the holding times $\{T_n:n\geq 1\}$ are i.i.d. with
\emph{generalized Mittag Leffler distribution}, i.e. with
continuous density $f_{\nu,h,\lambda}$ defined by
$$f_{\nu,h,\lambda}(t)=\lambda^h t^{\nu h-1}E_{\nu,\nu h}^h(-\lambda t^\nu)1_{(0,\infty)}(t).$$
We remark that, if we set $h=1$, we recover the same process in
\cite{BeghinOrsingher2009}-\cite{BeghinOrsingher2010} (see eq.
(2.16) in \cite{BeghinOrsingher2010}; see also eq. (4.14) in
\cite{BeghinOrsingher2009}). Moreover $f_{\nu,k,\lambda}$
coincides with eq. (2.19) in \cite{BeghinOrsingher2010}, where $k$
is integer. Finally we have
$f_{1,h,\lambda}(t)=\frac{\lambda^h}{\Gamma(h)}t^{h-1}e^{-\lambda
t}1_{(0,\infty)}(t)$ which is a Gamma density; thus we obtain the
classical case with exponentially distributed holding times for
$(\nu,h)=(1,1)$.

Now, in view of what follows, it is useful to recall that
\begin{equation}\label{eq:cumulant}
\kappa_{\nu,h,\lambda}(\theta):=\log\mathbb{E}[e^{\theta T_1}]=
\left\{\begin{array}{ll}
h\log\frac{\lambda}{\lambda+(-\theta)^\nu}&\ \mathrm{if}\ \theta\leq 0\\
\infty&\ \mathrm{if}\ \theta>0
\end{array}\right.\ \mathrm{for}\ \nu\in(0,1),
\end{equation}
and, for $\nu=1$,
$$\kappa_{1,h,\lambda}(\theta):=\log\mathbb{E}[e^{\theta T_1}]=
\left\{\begin{array}{ll}
h\log\frac{\lambda}{\lambda-\theta}&\ \mathrm{if}\ \theta<\lambda\\
\infty&\ \mathrm{if}\ \theta\geq \lambda.
\end{array}\right.$$

We conclude with the outline of this section. We start with the
LDPs for $\{\bar{T}_n:n\geq 1\}$, where
$\bar{T}_n:=\frac{T_1+\cdots+T_n}{n}$ for all $n\geq 1$, and for
$\left\{\frac{M_{\nu,h,\lambda}(t)}{t}:t>0\right\}$; moreover, for
the second LDP, we discuss the possible application of G\"{a}rtner
Ellis Theorem. In particular we study in detail the fractional
case $\nu=\frac{1}{2}$, for which we can provide explicit
expressions for the rate functions, and we recover the LDP for
$\left\{\frac{M_{\nu,1,\lambda}(t)}{t}:t>0\right\}$ (concerning
the case $h=1$) by taking into account that it can be expressed in
terms of a classical Poisson process computed at an independent
random time given by a reflecting Brownian motion with variance
parameter 2. Finally we present some results for the ruin
probabilities concerning an insurance model with a fractional
Poisson claim number process.

\subsection{The basic LDPs}\label{sub:basicLDPs}
We start with two LDPs which can be easily proved: the first one
(Proposition \ref{prop:LDP-holding-times}) concerns
$\{\bar{T}_n:n\geq 1\}$, and it is a particular case of Cram\'{e}r
Theorem (see e.g. Theorem 2.2.3 in \cite{DemboZeitouni}); the
second one (Proposition \ref{prop:LDP-counting-process}) concerns
$\left\{\frac{M_{\nu,h,\lambda}(t)}{t}:t>0\right\}$, and its proof
is based on the combination of the first result and known results
in the literature for nondecreasing processes and their inverses
(here we refer to \cite{DuffieldWhitt} which treats this kind of
problem in a wide generality, allowing non-linear scaling
functions; more precisely we refer to Theorem 1.1(i) in
\cite{DuffieldWhitt} with $u,v,w$ as the identity function because
$I_{\nu,h,\lambda}^{(T)}$ has no peaks with the unique base
$x=\infty$ if $\nu\in(0,1)$, and $x=\frac{h}{\lambda}$ if
$\nu=1$).

\begin{Proposition}\label{prop:LDP-holding-times}
The sequence $\{\bar{T}_n:n\geq 1\}$ satisfies the LDP with rate
function $I_{\nu,h,\lambda}^{(T)}$ defined by
$I_{\nu,h,\lambda}^{(T)}(x):=\sup_{\theta\in\mathbb{R}}\{\theta
x-\kappa_{\nu,h,\lambda}(\theta)\}$. In particular, for $\nu=1$,
we have
$$I_{1,h,\lambda}^{(T)}(x)=\left\{\begin{array}{ll}
h\left(\frac{\lambda x}{h}-1-\log\frac{\lambda x}{h}\right)&\ if\ x>0\\
\infty&\ if\ x\leq 0,
\end{array}\right.$$
and it is a good rate function. For $\nu\in(0,1)$ we have:
$I_{\nu,h,\lambda}^{(T)}(x)=\infty$ for $x\leq 0$,
$I_{\nu,h,\lambda}^{(T)}(x)$ is decreasing on $(0,\infty)$,
$\lim_{x\downarrow 0}I_{\nu,h,\lambda}^{(T)}(x)=\infty$,
$\lim_{x\to\infty}I_{\nu,h,\lambda}^{(T)}(x)=0$, the rate function
$I_{\nu,h,\lambda}^{(T)}$ is not good.
\end{Proposition}

We remark that, for all $\nu\in(0,1]$ and for all $h>0$, we have
$\kappa_{\nu,h,\lambda}(\theta)=h\kappa_{\nu,1,\lambda}(\theta)$
for all $\theta\in\mathbb{R}$, and therefore
$I_{\nu,h,\lambda}^{(T)}(x)=hI_{\nu,1,\lambda}^{(T)}(\frac{x}{h})$
for all $x\in\mathbb{R}$.

\begin{Proposition}\label{prop:LDP-counting-process}
The family $\left\{\frac{M_{\nu,h,\lambda}(t)}{t}:t>0\right\}$
satisfies the LDP with good rate function
$I_{\nu,h,\lambda}^{(M)}$ defined by
$$I_{\nu,h,\lambda}^{(M)}(x):=\left\{\begin{array}{ll}
xI_{\nu,h,\lambda}^{(T)}(1/x)&\ if\ x>0\\
\lambda 1_{\nu=1}&\ if\ x=0\\
\infty&\ if\ x<0.
\end{array}\right.$$
In particular, for $\nu=1$, we have
$$I_{1,h,\lambda}^{(M)}(x)=\left\{\begin{array}{ll}
hx\log\frac{hx}{\lambda}-hx+\lambda&\ if\ x\geq 0\\
\infty&\ if\ x<0.
\end{array}\right.$$
For $\nu\in(0,1)$ we have: $I_{\nu,h,\lambda}^{(M)}(x)=\infty$ for
$x<0$, $I_{\nu,h,\lambda}^{(M)}(x)$ is increasing on $[0,\infty)$,
$\lim_{x\downarrow
0}I_{\nu,h,\lambda}^{(M)}(x)=I_{\nu,h,\lambda}^{(M)}(0)=0$,
$\lim_{x\to\infty}I_{\nu,h,\lambda}^{(M)}(x)=\infty$.
\end{Proposition}

\paragraph{A discussion on G\"{a}rtner Ellis Theorem for the proof of Proposition \ref{prop:LDP-counting-process}.}
It is well-known that the rate function $I_{\nu,h,\lambda}^{(T)}$
is convex on $(0,\infty)$. Moreover $I_{\nu,h,\lambda}^{(M)}$ is
also convex on $(0,\infty)$; in fact, for $x_1,x_2\in(0,\infty)$
and $\gamma\in[0,1]$, we have
\begin{align*}
I_{\nu,h,\lambda}^{(M)}(\gamma x_1+(1-\gamma)x_2)=&(\gamma
x_1+(1-\gamma)x_2)I_{\nu,h,\lambda}^{(T)}\left(\frac{1}{\gamma
x_1+(1-\gamma)x_2}\right)\\
=&(\gamma
x_1+(1-\gamma)x_2)I_{\nu,h,\lambda}^{(T)}\left(\frac{\gamma
x_1}{\gamma x_1+(1-\gamma)x_2}\cdot\frac{1}{x_1}+\frac{(1-\gamma)
x_2}{\gamma x_1+(1-\gamma)x_2}\cdot\frac{1}{x_2}\right)
\end{align*}
and, by the convexity of $I_{\nu,h,\lambda}^{(T)}$, we get
\begin{align*}
I_{\nu,h,\lambda}^{(M)}(\gamma x_1+(1-\gamma)x_2)\leq&\gamma
x_1I_{\nu,h,\lambda}^{(T)}\left(\frac{1}{x_1}\right)+(1-\gamma)
x_2I_{\nu,h,\lambda}^{(T)}\left(\frac{1}{x_2}\right)\\
\leq&\gamma
I_{\nu,h,\lambda}^{(M)}(x_1)+(1-\gamma)I_{\nu,h,\lambda}^{(M)}(x_2).
\end{align*}
Then one could try to prove the LDP in Proposition
\ref{prop:LDP-counting-process} with an application of G\"{a}rtner
Ellis Theorem. If this was possible, we should have
$$\lim_{t\to\infty}\frac{1}{t}\log\mathbb{E}\left[e^{\theta
M_{\nu,h,\lambda}(t)}\right]=\left\{\begin{array}{ll}
\lambda(e^{\theta/h}-1)&\ \mathrm{if}\ \nu=1\\
(\lambda(e^{\theta/h}-1))^{1/\nu}1_{\theta\geq 0}&\ \mathrm{if}\
\nu\in(0,1)
\end{array}\right.=:\Lambda_{\nu,h,\lambda}(\theta)\ (\mbox{for all}\
\theta\in\mathbb{R})$$
and, since the function
$\theta\mapsto\Lambda_{\nu,h,\lambda}(\theta)$ satisfies the
hypotheses of G\"{a}rtner Ellis Theorem in both cases $\nu=1$ and
$\nu\in(0,1)$, we should get the LDP with rate function
$\Lambda_{\nu,h,\lambda}^*$ defined by
\begin{equation}\label{eq:rf-as-Legendre-transform}
\Lambda_{\nu,h,\lambda}^*(x)=\sup_{\theta\in\mathbb{R}}\left\{\theta
x-\Lambda_{\nu,h,\lambda}(\theta)\right\}\ (\mbox{for all}\
x\in\mathbb{R})
\end{equation}
because $\Lambda_{\nu,h,\lambda}^*$ coincides with the rate
function $I_{\nu,h,\lambda}^{(M)}$ in Proposition
\ref{prop:LDP-counting-process}. However we can have some
difficulties with this approach because we cannot have an
expression of the moment generating function
$\mathbb{E}\left[e^{\theta M_{\nu,h,\lambda}(t)}\right]$. For
completeness we also remark that, in both the cases $\nu=1$ and
$\nu\in(0,1)$, the function
$\theta\mapsto\Lambda_{\nu,h,\lambda}(\theta)$ above meets eq.
(12)-(13) in \cite{GlynnWhitt}, i.e.:
$$\Lambda_{\nu,h,\lambda}(\theta)=-\kappa_{\nu,h,\lambda}^{-1}(-\theta)\
\mathrm{for}\ \left\{\begin{array}{ll}
\theta\in\mathbb{R}&\ \mathrm{if}\ \nu=1\\
\theta\geq 0&\ \mathrm{if}\ \nu\in(0,1);
\end{array}\right.$$
$$\kappa_{\nu,h,\lambda}(\theta)=-\Lambda_{\nu,h,\lambda}^{-1}(-\theta)\
\mathrm{for}\ \left\{\begin{array}{ll}
\theta<\lambda&\ \mathrm{if}\ \nu=1\\
\theta\leq 0&\ \mathrm{if}\ \nu\in(0,1).
\end{array}\right.$$

\paragraph{Some remarks on the fractional case $\nu=\frac{1}{2}$.}
We can provide explicit formulas for the rate functions presented
above. By Proposition \ref{prop:LDP-holding-times}, we have
\begin{align*}
I_{\frac{1}{2},h,\lambda}^{(T)}(x)=&\sup_{\theta\leq
0}\left\{\theta
x-h\log\left(\frac{\lambda}{\lambda+(-\theta)^{\frac{1}{2}}}\right)\right\}=\left\{\theta
x-h\log\left(\frac{\lambda}{\lambda+(-\theta)^{\frac{1}{2}}}\right)
\right\}_{\theta=-\left(-\frac{\lambda}{2}+\frac{1}{2}\sqrt{\lambda^2+\frac{2h}{x}}\right)^2}\\
=&-\left(-\frac{\lambda}{2}+\frac{1}{2}\sqrt{\lambda^2+\frac{2h}{x}}\right)^2x-
h\log\left(\frac{\lambda}{\lambda-\frac{\lambda}{2}+\frac{1}{2}\sqrt{\lambda^2+\frac{2h}{x}}}\right)\\
=&-x\left(\frac{1}{2}\sqrt{\lambda^2+\frac{2h}{x}}-\frac{\lambda}{2}\right)^2+
h\log\left(\frac{1}{2}+\frac{1}{2}\sqrt{1+\frac{2h}{\lambda^2x}}\right)\
(\mbox{for all}\ x>0);
\end{align*}
thus, by Proposition \ref{prop:LDP-counting-process}, we have
\begin{equation}\label{eq:particular-fractional-case}
I_{\frac{1}{2},h,\lambda}^{(M)}(x)=xI_{\frac{1}{2},h,\lambda}^{(T)}(1/x)=hx\log\left(\frac{1}{2}+\frac{1}{2}\sqrt{1+\frac{2hx}{\lambda^2}}\right)-
\left(\frac{1}{2}\sqrt{\lambda^2+2hx}-\frac{\lambda}{2}\right)^2\
(\mbox{for all}\ x>0).
\end{equation}
We remark that $I_{\frac{1}{2},h,\lambda}^{(M)}$ in
\eqref{eq:particular-fractional-case} meets
$\Lambda_{\frac{1}{2},h,\lambda}^*$ in
\eqref{eq:rf-as-Legendre-transform} presented in the above
discussion on the proof of the LDP in Proposition
\ref{prop:LDP-counting-process} with an application of G\"{a}rtner
Ellis Theorem: the case $x\leq 0$ is trivial (the details are
omitted) and, for $x>0$, we have
\begin{align*}
\sup_{\theta\geq 0}\left\{\theta
x-(\lambda(e^{\theta/h}-1))^2\right\}=&\left\{\theta
x-(\lambda(e^{\theta/h}-1))^2
\right\}_{\theta=h\log\left(\frac{1}{2}+\frac{1}{2}\sqrt{1+\frac{2hx}{\lambda^2}}\right)}\\
=&h\log\left(\frac{1}{2}+\frac{1}{2}\sqrt{1+\frac{2hx}{\lambda^2}}\right)x-
\lambda^2\left(\frac{1}{2}\sqrt{1+\frac{2hx}{\lambda^2}}-\frac{1}{2}\right)^2\\
=&hx\log\left(\frac{1}{2}+\frac{1}{2}\sqrt{1+\frac{2hx}{\lambda^2}}\right)-
\left(\frac{1}{2}\sqrt{\lambda^2+2hx}-\frac{\lambda}{2}\right)^2=I_{\frac{1}{2},h,\lambda}^{(M)}(x).
\end{align*}

\paragraph{An alternative proof of Proposition \ref{prop:LDP-counting-process} for the case $(\nu,h)=(\frac{1}{2},1)$.}
The starting point consists of the following representation in the
literature (see Remark 2.1 in \cite{BeghinOrsingher2009}): for
each fixed $t>0$, $M_{\frac{1}{2},1,\lambda}(t)$ is distributed as
$N_\lambda(|B(2t)|)$ where
$$\left\{
\begin{array}{ll}
\{N_\lambda(t):t\geq 0\}\ \mbox{and}\ \{B(t):t\geq 0\}\ \mbox{are independent},\\
\{N_\lambda(t):t\geq 0\}\ \mbox{is a classical Poisson process, i.e. it is distributed as}\ \{M_{1,1,\lambda}(t):t\geq 0\},\\
\mbox{and}\ \{B(t):t\geq 0\}\ \mbox{is a standard Brownian motion}
\end{array}\right.$$
(note that $\{N_\lambda(|B(2t)|):t\geq 0\}$ does not represent a
version of $\{M_{\frac{1}{2},1,\lambda}(t):t\geq 0\}$ the process
$\{N_\lambda(|B(2t)|):t\geq 0\}$ is nondecreasing with respect to
$t$). We start with the following two statements.
\begin{enumerate}
\item The family of random variables $\left\{\frac{|B(2t)|}{t}:t>0\right\}$ satisfies
the LDP with good rate function $J$ defined by
$$J(y):=\left\{\begin{array}{ll}
\frac{y^2}{4}&\ \mathrm{if}\ y\geq 0\\
\infty&\ \mathrm{if}\ y<0.
\end{array}\right.$$
\emph{Sketch of the proof.} Firstly
$\left\{\frac{B(2t)}{t}:t>0\right\}$ satisfies the LDP by an easy
application of G\"{a}rtner Ellis Theorem with the good rate
function $H$ defined by $H(x)=\frac{x^2}{4}$; then the required
LDP holds by applying the contraction principle (see e.g. Theorem
4.2.1 in \cite{DemboZeitouni}) and noting that
$J(y)=\inf\{H(x):|x|=y\}$ for all $y\in\mathbb{R}$.
\item If $y_t\to y$ as $t\to\infty$, then
$\left\{\frac{N_\lambda(y_tt)}{t}:t>0\right\}$ satisfies the LDP
with rate function $K(\cdot|y)$ defined by
$$K(x|y):=\sup_{\theta\in\mathbb{R}}(\theta x-\lambda y(e^\theta-1))=\left\{
\begin{array}{ll}
\left\{\begin{array}{ll}
x\log\left(\frac{x}{\lambda y}\right)-x+\lambda y&\ \mathrm{if}\ x\geq 0\\
\infty&\ \mathrm{if}\ x<0,
\end{array}\right.&\ \mathrm{if}\ y>0,\\
\Delta_0(x)&\ \mathrm{if}\ y=0,
\end{array}\right.$$
where $\Delta_0$ is the function defined by
$$\Delta_0(x):=\left\{\begin{array}{ll}
0&\ \mathrm{if}\ x=0\\
\infty&\ \mathrm{if}\ x\neq 0.
\end{array}\right.$$
Moreover the function $(x,y)\mapsto K(x|y)$ is lower
semi-continuous on $[0,\infty)\times [0,\infty)$.\\
\emph{Sketch of the proof.} The LDP can be proved by an easy
application of G\"{a}rtner Ellis Theorem; moreover $(x,y)\mapsto
K(x|y)$ is a lower semi-continuous function on $[0,\infty)\times
[0,\infty)$ because, if $(x_n,y_n)\to (x,y)$, for all
$\theta\in\mathbb{R}$ we have
$$\liminf_{n\to\infty}K(x_n|y_n)\geq\liminf_{n\to\infty}(\theta x_n-\lambda y_n(e^\theta-1))=\theta x-\lambda y(e^\theta-1)$$
and we conclude by taking the supremum with respect to
$\theta\in\mathbb{R}$.
\end{enumerate}
In conclusion, by Theorem 2.3 in \cite{Chaganty} (namely we mean
the LDP for marginal distributions), the family of random
variables $\left\{\frac{N_\lambda(|B(2t)|)}{t}:t>0\right\}$
satisfies the LDP with rate function
$J_{\frac{1}{2},1,\lambda}^{(M)}$ (say) defined by
$$J_{\frac{1}{2},1,\lambda}^{(M)}(x):=\inf\{K(x|y)+J(y):y\geq 0\}.$$

Finally we show that $J_{\frac{1}{2},1,\lambda}^{(M)}$ coincides
with $I_{\frac{1}{2},1,\lambda}^{(M)}$. The equality
$J_{\frac{1}{2},1,\lambda}^{(M)}(x)=I_{\frac{1}{2},1,\lambda}^{(M)}(x)$
can be easily checked if $x<0$ (because $K(x|y)=\infty$ for all
$y\geq 0$) and if $x=0$ (because
$J_{\frac{1}{2},1,\lambda}^{(M)}(0)=\inf\left\{\lambda
y+\frac{y^2}{4}:y\geq 0\right\}=0$ since the infimum is attained
at $y=0$). If $x>0$ we have
$$J_{\frac{1}{2},1,\lambda}^{(M)}(x)=\inf\left\{x\log\left(\frac{x}{\lambda y}\right)-x+\lambda
y+\frac{y^2}{4}:y>0\right\};$$ then one can easily check that the
infimum is attained $y=\sqrt{\lambda^2+2x}-\lambda$, and we obtain
\begin{align*}
J_{\frac{1}{2},1,\lambda}^{(M)}(x)=&x\log\left(\frac{x}{\lambda(\sqrt{\lambda^2+2x}-\lambda)}\right)-x+
\lambda\left(\sqrt{\lambda^2+2x}-\lambda\right)+\frac{\left(\sqrt{\lambda^2+2x}-\lambda\right)^2}{4}\\
=&x\log\left(\frac{x(\sqrt{\lambda^2+2x}+\lambda)}{2\lambda
x}\right)-x+\lambda\left(\sqrt{\lambda^2+2x}-\lambda\right)+\left(\frac{1}{2}\sqrt{\lambda^2+2x}-\frac{\lambda}{2}\right)^2\\
=&x\log\left(\frac{1}{2}+\frac{1}{2}\sqrt{1+\frac{2x}{\lambda^2}}\right)-x
+\lambda\left(\sqrt{\lambda^2+2x}-\lambda\right)+\left(\frac{1}{2}\sqrt{\lambda^2+2x}-\frac{\lambda}{2}\right)^2;
\end{align*}
finally, by \eqref{eq:particular-fractional-case} with $h=1$, we
get
\begin{align*}
J_{\frac{1}{2},1,\lambda}^{(M)}(x)=&I_{\frac{1}{2},1,\lambda}^{(M)}(x)-x
+\lambda\left(\sqrt{\lambda^2+2x}-\lambda\right)+2\left(\frac{1}{2}\sqrt{\lambda^2+2x}-\frac{\lambda}{2}\right)^2\\
=&I_{\frac{1}{2},1,\lambda}^{(M)}(x)-x
+\lambda\sqrt{\lambda^2+2x}-\lambda^2+\frac{\lambda^2+2x}{2}+\frac{\lambda^2}{2}-\lambda\sqrt{\lambda^2+2x}
=I_{\frac{1}{2},1,\lambda}^{(M)}(x).
\end{align*}

\subsection{An insurance model with fractional Poisson claim number process}\label{sub:insurance}
In this subsection we study the ruin probability
$\Psi(u):=P(\{\exists t\geq 0:R(t)<0\})$ concerning the insurance
model
$$R(t):=u+ct-\sum_{k=1}^{M_{\nu,h,\lambda}(t)}U_k,$$
where (we refer to the terminology for eq. (5.1.14) in
\cite{RSST}) $\{R(t):t\geq 0\}$ is the reserve process, $u>0$ is
the initial capital of the company, $c>0$ is the premium rate and
$\{U_k:k\geq 1\}$ are the claim sizes assumed to be i.i.d.
positive random variables and independent of the claim number
process $\{M_{\nu,h,\lambda}(t):t\geq 0\}$ defined by
\eqref{eq:def-main-version}. Here we consider a slightly different
notation for the holding times, which will be denoted by
$\{T_n^{(\nu)}:n\geq 1\}$ instead of $\{T_n:n\geq 1\}$.

We always consider a fractional claim number process, i.e. we
assume that $\nu\in(0,1)$ and $h>0$. We recall that, if
$(\nu,h)=(1,1)$, the claim number process is a homogeneous Poisson
process and we have the \emph{compound Poisson model} (see e.g.
section 5.3 in \cite{RSST}; see also the \emph{Cram\'{e}r-Lundberg
model} in section 1.1 in \cite{EKM}).

It is easy to check that the ruin probability $\Psi(u)$ coincides
with a level crossing probability for the random walk
$\{\sum_{k=1}^n(U_k-cT_k^{(\nu)}):n\geq 1\}$, i.e.
$$\Psi(u)=P\left(\left\{\exists n\geq 1:\sum_{k=1}^n(U_k-cT_k^{(\nu)})>u\right\}\right);$$
this happens because the ruin can occur only at the time epochs of
the claims. Furthermore it is known that, if we consider the case
$\nu=1$, the ruin problem is non-trivial (i.e. $\Psi(u)\in(0,1)$)
if $c$ is large enough to have $\mathbb{E}[U_1-cT_1^{(1)}]<0$,
i.e. if the \emph{net profit condition}
$c>\frac{\lambda}{h}\mathbb{E}[U_1]$ holds (note that, for $h=1$,
this meets eq. (5.3.2) in \cite{RSST}, or eq. (1.7) in \cite{EKM}
(page 26)). On the contrary, for the fractional case $\nu\in(0,1)$
considered here, the ruin problem is non trivial for any $c>0$
because we have $\mathbb{E}[U_1-cT_1^{(\nu)}]=-\infty$.

Here we present two results which can be derived from
straightforward applications of Theorems 1-2 in
\cite{LehtonenNyrhinen} for the random walk
$\{\sum_{k=1}^n(U_k-cT_k^{(\nu)}):n\geq 1\}$, respectively. Thus
we need to consider the function $\tilde{\kappa}_\nu$ defined by
\begin{equation}\label{eq:global-cumulant}
\tilde{\kappa}_\nu(\theta):=\log\mathbb{E}[e^{\theta
U_1}]+\log\mathbb{E}[e^{-c\theta T_1^{(\nu)}}]
=\log\mathbb{E}[e^{\theta U_1}]+\kappa_{\nu,h,\lambda}(-c\theta),
\end{equation}
and the following condition:\\
$\mathbf{(C1)}$: \emph{there exists}
$w_{\nu,h,\lambda}\in(0,\infty)\cap\{\theta\in\mathbb{R}:\tilde{\kappa}_\nu(\theta)<\infty\}^\circ$
\emph{such that} $\tilde{\kappa}_\nu(w_{\nu,h,\lambda})=0$.\\

We start with the first result which provides an asymptotic
estimate of $\Psi(u)$ in the fashion of large deviations.

\begin{Proposition}\label{prop:LN1}
Assume that $\mathbf{(C1)}$ holds. Then we have
$\lim_{u\to\infty}\frac{1}{u}\log\Psi(u)=-w_{\nu,h,\lambda}$.
\end{Proposition}

Note that, if $\nu_1<\nu_2$ (with $\nu_1,\nu_2\in(0,1)$), then
$w_{\nu_1,h,\lambda}<w_{\nu_2,h,\lambda}$; this can be checked
noting that, by \eqref{eq:global-cumulant} and the definition of
$\kappa_{\nu,h,\lambda}$ in \eqref{eq:cumulant},
$\tilde{\kappa}_{\nu_1}(\theta)>\tilde{\kappa}_{\nu_2}(\theta)$
for $\theta>0$. Thus the smaller is the value $\nu$, the more
dangerous is the situation (i.e. the more slowly the ruin
probabilities decay as $u\to\infty$).\\

The second result gives an optimal importance sampling
distribution for the estimation of $\Psi(u)$ by simulation for
large values of $u$. We need some further preliminaries. Let
$P_U\otimes P_T$ be the common law for the random variables
$\{(U_n,T_n^{(\nu)}):n\geq 1\}$. Moreover, for each $\theta$ such
that $\tilde{\kappa}_\nu(\theta)<\infty$, let $P_U^\theta\otimes
P_T^\theta$ be the absolutely continuous law with density
$\frac{d(P_U^\theta\otimes P_T^\theta)}{d(P_U\otimes
P_T)}(x,t)=\frac{dP_U^\theta}{dP_U}(x)\frac{dP_T^\theta}{dP_T}(t)$
where $\frac{dP_U^\theta}{dP_U}(x)=e^{\theta
x-\log\mathbb{E}[e^{\theta U_1}]}$ and
$\frac{dP_T^\theta}{dP_T}(t)=e^{-c\theta
t-\log\mathbb{E}[e^{-c\theta T_1^{(\nu)}}]}=e^{-c\theta
t-\kappa_{\nu,h,\lambda}(-c\theta)}$. Here we think to have i.i.d.
random variables $\{(U_n,T_n^{(\nu)}):n\geq 1\}$ whose common law
is $P_U^\theta\otimes P_T^\theta$ (for some $\theta$); thus, in
particular, each one of the random variables
$\{(U_n,T_n^{(\nu)}):n\geq 1\}$ has independent components as
happens under the original law $P_U\otimes P_T$ (i.e.
$P_U^0\otimes P_T^0$) of the random variables.

\begin{Proposition}\label{prop:LN2}
Assume that $\mathbf{(C1)}$ holds. Then, for
$\theta=w_{\nu,h,\lambda}$, $P_U^\theta\otimes P_T^\theta$ is an
optimal importance sampling distribution for the estimation of
$\Psi(u)$ by simulation for large values of $u$.
\end{Proposition}

Note that the exponential change of measure $P_U^\theta\otimes
P_T^\theta$ presented above can be considered also for $\nu=1$.
Then we have the two following situations.
\begin{itemize}
\item If $\nu=1$, for $\theta>-\frac{\lambda}{c}$ we have
$$dP_T^\theta(t)=\frac{e^{-c\theta t}\frac{\lambda^h}{\Gamma(h)}t^{h-1}e^{-\lambda t}1_{(0,\infty)}(t)dt}
{\int_0^\infty e^{-c\theta
y}\frac{\lambda^h}{\Gamma(h)}y^{h-1}e^{-\lambda y}dy}=
\frac{(c\theta+\lambda)^h}{\Gamma(h)}t^{h-1}e^{-(c\theta+\lambda)t}1_{(0,\infty)}(t)dt;$$
thus $\{P_T^\theta:\theta>-\frac{\lambda}{c}\}$ are all Gamma
distributions.
\item If $\nu\in(0,1)$, for $\theta\geq 0$ we have
$$dP_T^\theta(t)=\frac{e^{-c\theta t}\lambda^h t^{\nu h-1}E_{\nu,\nu h}^h(-\lambda t^\nu)1_{(0,\infty)}(t)dt}
{\int_0^\infty e^{-c\theta y}\lambda^h y^{\nu h-1}E_{\nu,\nu
h}^h(-\lambda y^\nu)dy}=e^{-c\theta
t}(\lambda+(c\theta)^\nu)^ht^{\nu h-1}E_{\nu,\nu h}^h(-\lambda
t^\nu)1_{(0,\infty)}(t)dt;$$ thus $\{P_T^\theta:\theta>0\}$ are
not generalized Mittag Leffler distributions as it is $P_T^0$
because the equality $e^{-c\theta t}E_{\nu,\nu h}^h(-\lambda
t^\nu)=E_{\nu,\nu h}^h(-(\lambda+(c\theta)^\nu)t^\nu)$ holds if
and only if $\theta=0$ (on the contrary the equality always holds
if $\nu=1$).
\end{itemize}

\section{Results for the alternative version (weighted Poisson laws)}\label{sec:alternative-version}
In this section we consider an alternative version of the
fractional Poisson process $\{A_{\nu,\lambda}(t):t\geq 0\}$ which
is the first version presented in section 4 in
\cite{BeghinOrsingher2009} with $t^\nu$ in place of $t$:
$$P(A_{\nu,\lambda}(t)=k)=\frac{(\lambda t^\nu)^k}{\Gamma(\nu
k+1)}\frac{1}{E_{\nu,1}(\lambda t^\nu)}\ \mbox{for all}\
k\in\mathbb{N}^*:=\{0,1,2,3,\ldots\}.$$ We remark that each random
variable $A_{\nu,\lambda}(t)$ has a particular weighted Poisson
distribution (we refer to the terminology in
\cite{JohnsonKotzKemp}, page 90; see also the references cited
therein), and the weights do not depend on $t$. More precisely,
for each fixed $t$, the discrete density of $A_{\nu,\lambda}(t)$
is
$$q_w(k):=\frac{w(k)q(k)}{\sum_{j\geq 0}w(j)q(j)}\ \mbox{for all}\ k\in\mathbb{N}^*,$$
where the density $\{q(k):k\in\mathbb{N}^*\}$ and the weights
$\{w(k):k\in\mathbb{N}^*\}$ are defined by $q(k):=\frac{(\lambda
t^\nu)^k}{k!}e^{-\lambda t^\nu}$ (the classical Poisson density
with mean $\lambda t^\nu$) and $w(k):=\frac{k!}{\Gamma(\nu k+1)}$,
respectively.

In this section we prove the  LDP for
$\left\{\frac{A_{\nu,\lambda}(t)}{t}:t>0\right\}$ and we provide a
formula (see eq. \eqref{eq:*} below) for the rate function in
terms of a suitable limit of normalized relative entropies (see
eq. \eqref{eq:limit-normalized-relative-entropies} below). This
has some analogy with a recent result for stationary Gaussian
processes (see section 2 in \cite{MacciPetrella}); moreover it is
well-known (see e.g. the discussion in \cite{Varadhan}) that the
rate functions are often expressed in terms of the relative
entropy.

\begin{Proposition}\label{prop:LDP-alternative}
For $\nu\in(0,1]$,
$\left\{\frac{A_{\nu,\lambda}(t)}{t}:t>0\right\}$ satisfies the
LDP with good rate function $I_{\nu,\lambda}^{(A)}$ defined by
$$I_{\nu,\lambda}^{(A)}(x):=\left\{\begin{array}{ll}
\nu x\log\frac{\nu x}{\lambda^{1/\nu}}-\nu x+\lambda^{1/\nu}&\ if\ x\geq 0\\
\infty&\ if\ x<0.
\end{array}\right.$$
\textbf{Remark.} For each fixed $t\geq 0$, $A_{1,\lambda}(t)$ is
distributed as $M_{1,1,\lambda}(t)$ in
\eqref{eq:def-main-version}. Thus, if $\nu=1$, we recover
Proposition \ref{prop:LDP-counting-process} with $h=1$ and,
actually, one can check that $I_{1,\lambda}^{(A)}$ coincides with
$I_{1,1,\lambda}^{(M)}$.
\end{Proposition}
\noindent\emph{Proof.} Firstly we can immediately check that
$$\mathbb{E}\left[e^{\theta
A_{\nu,\lambda}(t)}\right]=\frac{E_{\nu,1}(e^\theta\lambda
t^\nu)}{E_{\nu,1}(\lambda t^\nu)}$$ for all $\theta\in\mathbb{R}$;
note that $\mathbb{E}\left[e^{\theta
A_{\nu,\lambda}(t)}\right]=m(e^\theta)$, where $m(\cdot)$ is the
probability generating function in eq. (4.4) in
\cite{BeghinOrsingher2009} (with $t^\nu$ in place of $t$).
Therefore, by using \eqref{eq:asymptotic}, we can check the limit
$$\lim_{t\to\infty}\frac{1}{t}\log\mathbb{E}\left[e^{\theta A_{\nu,\lambda}(t)}\right]=\lambda^{1/\nu}(e^{\theta/\nu}-1).$$
Then, by G\"{a}rtner Ellis Theorem, the LDP holds with good rate
function $I_{\nu,\lambda}^{(A)}$ defined by
$I_{\nu,\lambda}^{(A)}(x):=\sup_{\theta\in\mathbb{R}}\{\theta
x-\lambda^{1/\nu}(e^{\theta/\nu}-1)\}$ which coincides with the
rate function in the statement (we omit the details). $\Box$\\

In view of what follows we recall the definition and some
properties of the relative entropy (see e.g. section 2.3 in
\cite{CoverThomas}). Given two probability measures $Q_1$ and
$Q_2$ on the same measurable space $(\Omega,\mathcal{B}(\Omega))$,
we write $Q_1\ll Q_2$ to mean that $Q_1$ is absolutely continuous
with respect to $Q_2$ and, in such a case, the density will be
denoted by $\frac{dQ_1}{dQ_2}$. Then the relative entropy of $Q_1$
with respect to $Q_2$ is defined by
$$H(Q_1|Q_2)=\left\{\begin{array}{ll}
\int_\Omega\log(\frac{dQ_1}{dQ_2}(\omega))Q_1(d\omega)&\ \mathrm{if}\ Q_1\ll Q_2\\
\infty &\ \mathrm{otherwise}.
\end{array}\right.$$
It is known that $H(Q_1|Q_2)$ is nonnegative and it is equal to
zero if and only if $Q_1=Q_2$.

Now, in view of what follows, let $Q_{\nu,\lambda,t}$ be the law
of $A_{\nu,\lambda}(t)$; here we also allow the case $\lambda=0$,
and $Q_{\nu,0,t}$ is the law of the constant random variable equal
to 0. Then, if we consider the following limit of normalized
relative entropies
\begin{equation}\label{eq:limit-normalized-relative-entropies}
\mathcal{H}_\nu(\lambda_1|\lambda_2):=\lim_{t\to\infty}\frac{1}{t}H(Q_{\nu,\lambda_1,t}|Q_{\nu,\lambda_2,t})
\end{equation}
(for $\nu\in(0,1]$ and $\lambda_1,\lambda_2\geq 0$), we have
\begin{equation}\label{eq:*}
I_{\nu,\lambda}^{(A)}(x)=\mathcal{H}_\nu((\nu x)^\nu|\lambda)\
\mbox{for all}\ x\geq 0
\end{equation}
as an immediate consequence of the following result.

\begin{Proposition}\label{prop:formula}
For $\nu\in(0,1]$ and $\lambda_1,\lambda_2\geq 0$, we have
$\mathcal{H}_\nu(\lambda_1|\lambda_2)=\lambda_1^{1/\nu}\log\frac{\lambda_1^{1/\nu}}{\lambda_2^{1/\nu}}-\lambda_1^{1/\nu}+\lambda_2^{1/\nu}$.
\end{Proposition}
\noindent\emph{Proof.} We start assuming that
$\lambda_1,\lambda_2>0$. We have the following chain of equalities
where, for the latter equality, we take into account eq. (4.6) in
\cite{BeghinOrsingher2009} (with $t^\nu$ in place of $t$) for the
expected value $\sum_{k=0}^\infty kP(A_{\nu,\lambda_1}(t)=k)$:
\begin{align*}
\frac{1}{t}H(Q_{\nu,\lambda_1,t}|Q_{\nu,\lambda_2,t})=&\frac{1}{t}\sum_{k=0}^\infty
P(A_{\nu,\lambda_1}(t)=k)\log\left(\frac{P(A_{\nu,\lambda_1}(t)=k)}{P(A_{\nu,\lambda_2}(t)=k)}\right)\\
=&\frac{1}{t}\sum_{k=0}^\infty
P(A_{\nu,\lambda_1}(t)=k)\log\left(\frac{\lambda_1^k}{\lambda_2^k}
\frac{E_{\nu,1}(\lambda_2 t^\nu)}{E_{\nu,1}(\lambda_1
t^\nu)}\right)\\
=&\frac{1}{t}\log\frac{\lambda_1}{\lambda_2}\sum_{k=0}^\infty
kP(A_{\nu,\lambda_1}(t)=k)+\frac{1}{t}\log\left(\frac{E_{\nu,1}(\lambda_2
t^\nu)}{E_{\nu,1}(\lambda_1 t^\nu)}\right)\\
=&\frac{1}{t}\frac{\lambda_1t^\nu}{\nu}\frac{E_{\nu,\nu}(\lambda_1t^\nu)}{E_{\nu,1}(\lambda_1t^\nu)}\cdot\log\frac{\lambda_1}{\lambda_2}
+\frac{1}{t}\log\left(\frac{E_{\nu,1}(\lambda_2t^\nu)}{E_{\nu,1}(\lambda_1t^\nu)}\right).
\end{align*}
Then, by using \eqref{eq:asymptotic}, the limit in
\eqref{eq:limit-normalized-relative-entropies} holds with
$$\mathcal{H}_\nu(\lambda_1|\lambda_2)=
\frac{\lambda_1^{1/\nu}}{\nu}\log\frac{\lambda_1}{\lambda_2}+\lambda_2^{1/\nu}-\lambda_1^{1/\nu}
=\lambda_1^{1/\nu}\log\frac{\lambda_1^{1/\nu}}{\lambda_2^{1/\nu}}-\lambda_1^{1/\nu}+\lambda_2^{1/\nu}.$$
Thus the proof of the proposition is complete when
$\lambda_1,\lambda_2>0$, and now we give some details for the
other cases. If $\lambda_1=0$ and $\lambda_2>0$, we can consider
this procedure, but the above sum reduces to the first addendum
(the one with $k=0$) and we have
$\mathcal{H}_\nu(\lambda_1|\lambda_2)=\lambda_2^{1/\nu}$. If
$\lambda_2=0$, we have
$$\mathcal{H}_\nu(\lambda_1|0)=\left\{\begin{array}{ll}
0&\ \mathrm{if}\ \lambda_1=0\\
\infty&\ \mathrm{if}\ \lambda_1>0
\end{array}\right.$$
because, for all $t>0$, we trivially have
$H(Q_{\nu,0,t}|Q_{\nu,0,t})=0$ and, if $\lambda_1>0$,
$H(Q_{\nu,\lambda_1,t}|Q_{\nu,0,t})=\infty$. $\Box$\\

Finally we remark that
\begin{align*}
\frac{1}{t}H(Q_{1,\lambda_1,t}|Q_{1,\lambda_2,t})&=\frac{1}{t}\sum_{k=0}^\infty
P(A_{1,\lambda_1}(t)=k)\log\left(\frac{\lambda_1^k}{\lambda_2^k}
\frac{E_{1,1}(\lambda_2 t)}{E_{1,1}(\lambda_1 t)}\right)\\
&=\frac{1}{t}\log\frac{\lambda_1}{\lambda_2}\sum_{k=0}^\infty
k\frac{(\lambda_1t)^k}{k!}e^{-\lambda_1t}+\frac{1}{t}\log\left(e^{(\lambda_2-\lambda_1)t}\right)\\
&=\lambda_1\log\frac{\lambda_1}{\lambda_2}-\lambda_1+\lambda_2=H(Q_{1,\lambda_1,1}|Q_{1,\lambda_2,1})
\end{align*}
does not depend on $t>0$, and therefore coincides with
$\mathcal{H}_1(\lambda_1|\lambda_2)$.

\end{document}